\documentclass[11pt, a4paper]{amsart}

\usepackage[english]{babel}
\usepackage{amsmath,enumerate, amsfonts, amssymb,amsthm, xy}
\usepackage[dvips]{graphics} 

\newtheorem{thm}{Theorem}[section]
\newtheorem{lemma}[thm]{Lemma}
\newtheorem{prop}[thm]{Proposition}
\newtheorem{cor}[thm]{Corollary}

\theoremstyle{definition}
\newtheorem{defi}[thm]{Definition}

\newtheorem{ex}[thm]{Example}
\newtheorem{rmk}[thm]{Remark}

\DeclareMathOperator{\ord}{ord}
\DeclareMathOperator{\jac}{jac}

\DeclareMathOperator{\point}{point}

\title{Blow-Nash types of simple singularities}

\author{Goulwen Fichou}
\address {Institut Math\'ematiques de Rennes, Universit\'e de Rennes 1,
Campus de Beaulieu, 35042 Rennes Cedex, France}
\email{goulwen.fichou@univ-rennes1.fr}

\subjclass[2000]{Primary 14B05; Secondary 14P20, 14P25, 32S15.}
\keywords{Blow-Nash equivalence - Simple singularities - Virtual Poincar\'e polynomial}

\begin{document}

\begin{abstract} We address the question of the classification under blow-Nash equivalence of simple Nash function germs. We state that this classification coincides with the real analytic classification. We prove moreover that a simple germ can not be blow-Nash equivalent to a nonsimple one. The method is based on the computation of relevant coefficients of the real zeta functions associated to a Nash germ via motivic integration.
\end{abstract}

\maketitle
%%%%%%%%%%%%%%%%%%%%%%%%%%%%%%%%%%%%%%%%%%%%%%%%%%%%%%%%%%%%%%%%%%%%%%%%%%%%%%%%%%%%%%%%%%%%%%%%%%%%%%%%%%%%%%
The classification of real analytic singularities is a fascinating topic. In particular the choice of a good equivalence relation to study is a crucial point. In this paper we are concerned with the classification of real analytic and semi-algebraic function germs, called Nash function germs \cite{Shiota}, under blow-Nash equivalence \cite{fichou}. Nash function germs are blow-Nash equivalent if they become real analytic and semi-algebraically equivalent after some resolution of their singularities (see Definition \ref{defibN} for a precise statement). The study of such an equivalence relation, typical to the real setting, is motivated by the work of T.C. Kuo \cite{Kuo}. He noticed that, if the Whitney family $w(x,y)=xy(y-x)(y-tx)$ of real function germs at the origin admits infinitely many different analytic types for $t\in (0,1)$, it becomes analytically trivial after blowing-up the origin. He introduced in that way the notion of blow-analytic equivalence whose semi-algebraic counterpart is the blow-Nash equivalence. For more information concerning blow-analytic and blow-Nash equivalences we refer to the survey \cite{FP}.

The classification of complex and real analytic function germs have already a rich history and many invariants are associated with such germs. The classification begins with regarding at the corank, and index in the real context, of the second differential. Then we distinguish between more or less basic singularity types. In the classical book \cite{Arnold} we find notably the classification of all possible types of singularity up to modality two, with some complete lists given with normal forms of the singularities. Among them, the most simple singularities of modality zero (who admit only finitely many different analytic types in a sufficiently small neighbourhood) are the well-known $ADE$-singularities. Note that this classification still holds with Nash function germs by Nash Approximation Theorem (\cite{Shiota} and section \ref{ADE}).

In \cite{fichou} we have already established the invariance of corank and index under blow-Nash equivalence. This result was the first step toward a classification of Nash function germs with respect to the blow-Nash equivalence. The analog result concerning blow-analytic equivalence is still open.

In this paper we go further in the classification and consider simple Nash function germs. We know that, in general, the classification of real analytic (respectively Nash) function germs does not coincide under analytic and blow-analytic equivalence (respectively Nash or blow-Nash equivalence) as illustrated by Whitney example. Nevertheless, it is reasonable to hope that, for sufficiently simple singularities, these classifications may coincide. Theorem \ref{thm1} states that this is actually the case for simple Nash germs under blow-Nash equivalence. Namely, simple Nash germs, i.e. Nash germs belonging to the list of $ADE$-singularities, are blow-Nash equivalent if and only if they are analytically equivalent.

Actually, we prove even more. The blow-Nash classification not only coincides with the analytic classification for simple Nash germs, but moreover Theorem \ref{thm2} asserts that a nonsimple Nash germ can not belong to the same blow-Nash class as that of a simple germ!

We knew already that blow-Nash equivalence had a nice behaviour with respect to moduli. In particular a Nash family with isolated singularities admits locally only finitely many different blow-Nash types (\cite{fichou}, see also \cite{FP1,Koike}). Theorems \ref{thm1} and \ref{thm2} claim that this equivalence relation has also a very nice behaviour in terms of classification results.

The paper is organised as follows. In section 1 we recall the definition of the zeta functions associated with a Nash function germ \cite{fichou}. These invariants of blow-Nash equivalence, coming from the theory of motivic integration \cite{DL}, are formal power series whose coefficients are polynomials constructed by considering some measure over spaces of arcs connected to the given germ. The proof of the invariance of the index of corank under blow-Nash equivalence in \cite{fichou2} consisted in analysing the informations inclosed in the second coefficient of these formal power series. In this paper we need to inspect deeply inside the zeta functions, the informations coding the belonging to a simple singularity type (section \ref{ADE}). We put particular emphasis in section \ref{qua} on the case of quadratic polynomials that will be crucial in the remaining of the paper. Sections \ref{sectionsimple} and \ref{nonsimple} are devoted to the statement and the proof of the main theorems.

{\bf Acknowledgments.} The author is greatly indebted to T. Fukui for motivating discussions during his nice stay at Saitama University where this paper has been written. He is also grateful to the Japan Society for the Promotion of Science for its financial support.
%%%%%%%%%%%%%%%%%%%%%%%%%%%%%%%%%%%%%%%%%%%%%%%%%%%%%%%%%%%%%%%%%%%%%%%%%%%%%%%%%%%%%%%%%%%%%%%%%%%%%%%%%
\section{Blow-Nash equivalence}
The blow-Nash equivalence between Nash function germs is the analog in the Nash setting of the blow-analytic equivalence of T.C. Kuo \cite{Kuo}. Note that a Nash function (map, isomorphism, space) means a semi-algebraic and real analytic function (map, isomorphism, space). Blow-analytic equivalence has been proven to be a relevant equivalence relation between real analytic function germs (see \cite{FP} for a recent survey). Blow-Nash equivalence has been introduced in \cite{fichou} as a more algebraic version that enables us to use algebraic geometry, such as motivic integration \cite{DL}, for its study.

We first state its definition and then recall some important properties.

\begin{defi}\label{defibN}\begin{flushleft}\end{flushleft}
\begin{enumerate}
\item A Nash function germ $f:(\mathbb R^d,0) \longrightarrow (\mathbb R,0)$ admits $\pi:(M,\pi^{-1}(0))\longrightarrow (\mathbb R^d,0)$ as a Nash modification if:
\begin{enumerate}
\item[-] $\pi$ is a proper surjective Nash map between semi-algebraic neighbourhoods of $\pi^{-1}(0)$ in a Nash space $M$ and of $0$ in $\mathbb R^d$,
\item[-] $\pi$ is an isomorphism over the complement of the zero locus of $f$,
\item[-] the complexification of $\pi$ is an analytic isomorphism except on some thin subset,
\item[-] $f\circ \pi$ and $\jac \pi$ have only normal crossings simultaneously.
\end{enumerate}
\item Nash function germs $f,g:(\mathbb R^d,0) \longrightarrow (\mathbb R,0)$ are blow-Nash equivalent if there exists a blow-Nash isomorphism $h$ of $(\mathbb R^d,0)$ such that $f=g\circ h$, which means $h$ is a homeomorphism and moreover there exist Nash modifications $\pi_f:(M_f,E_f)\longrightarrow (\mathbb R^d,0)$ for $f$ and $\pi_g:(M_g,E_g)\longrightarrow (\mathbb R^d,0)$ for $g$ and a Nash isomorphism $H$ between $(M_f,E_f)$ and $(M_g,E_g)$, which is an isomorphism between the critical loci $E_f$ and $E_g$ considered as Nash spaces, such that $h\circ \pi_f=\pi_g\circ H$.
\end{enumerate}
\end{defi}

Similarly to the case of the blow-analytic equivalence, to find the right definition of the blow-Nash equivalence is still a work in progress. In particular this definition, where we impose a strong control on the critical loci, is closer to that of \cite{FKP} than to the original one of \cite{Kuo}. Note that despite its technical aspect, the blow-Nash equivalence has  nice geometric properties. For instance there are no moduli for Nash family with isolated singularity \cite{fichou}. Moreover the singular loci of the function germs are preserved \cite{P}. In particular, a function germ with an isolated singularity can not be equivalent to a function germ with a nonisolated singularity (this fact will be useful in section \ref{nonsimple}).

The behaviour of the blow-Nash equivalence seems to be better understood by looking at arcs that cross the singularity at the origin. In that spirit, the first invariant constructed for such an equivalence relation was the set of possible orders of series obtained by composition of a given function germ with analytic arcs passing through the origin. This set, known as Fukui Invariant \cite{IKK}, admits a natural generalisation via motivic integration \cite{DL}. More precisely not only the set of orders of arcs are invariant, but also some measure of the space of arcs that realises a given order \cite{fichou}. The associated invariants, called zeta functions, will play a crucial role in this paper.

We recall now the construction of these zeta functions. First we recall the definition and basic properties of the virtual Poincar\'e polynomial of a constructible real algebraic set that acts as the measure mentioned before.

\begin{prop}(\cite{MCP,fichou})\label{beta-nash} Take $i\in \mathbb N \cup 0$. The Betti number $b_i(\cdot)=\dim H_i(\cdot, \frac {\mathbb Z}{2 \mathbb Z})$, considered on compact nonsingular real algebraic sets, admits an unique extension as an additive map $\beta_i$ to the category of constructible real algebraic sets, with values in $\mathbb Z$. Namely
$$\beta(X)=\beta(Y)+\beta(X \setminus Y)$$
for $Y \subset X$ a closed subvariety of $X$. The polynomial $\beta(\cdot)=\sum_{i \geq 0} \beta_i(\cdot)u^i$ with values in $\mathbb Z [u]$ is multiplicative: $\beta(X\times Y)=\beta(X)\beta(Y)$ for constructible sets $X,Y$.
Finally $\beta(X)=\beta(Y)$ for Nash isomorphic constructible real algebraic sets $X$ and $Y$.
\end{prop}

The invariant $\beta_i$ is called the i-$th$ virtual Betti number, and the
polynomial $\beta$ the virtual Poincar\'e polynomial. By evaluation of the virtual Poincar\'e polynomial at $-1$ one recovers the Euler characteristic with compact support \cite{MCP}.

The following simple example illustrates the way to compute, in actual practice, the virtual Poincar\'e polynomial. We will perform more involved computations using intensively the additivity property in sections \ref{qua} and \ref{ADE}.

\begin{ex} If $\mathbb P^k$ denotes the real projective space of dimension
  $k$, which is nonsingular and compact, then $\beta(\mathbb
  P^k)=1+u+\cdots+u^k$ since $\dim H_i(\mathbb
  P^k, \frac {\mathbb Z}{2 \mathbb Z})=1$ for $i \in \{0,\ldots,k\}$ and $\dim H_i(\mathbb
  P^k, \frac {\mathbb Z}{2 \mathbb Z})=0$ otherwise. Now, compactify the affine line $\mathbb  A_{\mathbb R}^1$ in $\mathbb P^1$
  by adding one point at the infinity. By additivity $\beta(\mathbb  A_{\mathbb R}^1)=\beta(\mathbb P^1)-\beta(\point)=u,$ and so $\beta(\mathbb  A_{\mathbb R}^k)=u^k$ by multiplicativity.
\end{ex}

Now we are in position to introduce the zeta functions associated with a Nash function germ $f:(\mathbb R ^d,0) \longrightarrow (\mathbb R,0)$. Let $n \in \mathbb N$. Denote by $P_n[t]$ the set of polynomial arcs with coefficients in $\mathbb R^d$ of degree at most $n$ vanishing at the origin. We construct the naive zeta function and the zeta functions with sign of $f$ by considering the measure under the virtual Poincar\'e polynomial of the spaces of polynomial arcs $A_n(f),A_n^+(f)$ and $A_n^-(f)$ defined by:
$$A_n(f) =\{\gamma(t) \in P_n[t]: \ord_t f\circ \gamma(t) =n \},$$
$$A_n^{\pm1}(f)=\{\gamma(t) \in P_n[t]:f\circ \gamma (t)=\pm t^n+\cdots\}\subset A_n(f)$$
where the dots mean higher order terms.

The naive zeta function $Z_f(T)$ of $f$ is the formal power series in $\mathbb Z[u][[T]]$ defined by: 
$$Z_f(T)= \sum _{n \geq 1}{\beta (A_n(f))T^n},$$
whereas the zeta functions with sign are defined by:
$$Z_f^{+1}(T)= \sum _{n \geq 1}{\beta (A_n^{+1}(f))T^n}
\textrm{~~and~~ }Z_f^{-1}(T)= \sum _{n \geq 1}{\beta
(A_n^{-1}(f))T^n}.$$

\begin{thm}(\cite{fichou})\label{inv} Blow-Nash equivalent Nash function
germs share the same naive zeta function and the same zeta functions with sign. 
\end{thm}

\begin{rmk} The original definition of the zeta functions involves a correcting term of the form $u^{-nd}$ in front of $T^n$ which is important in the proof of Theorem \ref{inv} \cite{fichou}. We omit it in this paper for simplicity.
\end{rmk}

In the previous paper \cite{fichou2}, we used the invariance under blow-Nash equivalence of the $T^2$-coefficient of the zeta functions to prove the first step in the classification of Nash function germs under blow-Nash equivalence:

\begin{cor}(\cite{fichou2})\label{invcorank} The corank and index of a Nash function germ are invariant under blow-Nash equivalence.
\end{cor}
%%%%%%%%%%%%%%%%%%%%%%%%%%%%%%%%%%%%%%%%%%%%%%%%%%%%%%%%%%%%%%%%%%%%%%%%%%%%%%%%%%%%%%%%%%%%%%%%%%%%%%%%%%%%%%

\section{Quadratic case}\label{qua}
This section is devoted to the computation of the virtual Poincar\'e polynomial of some spaces of arcs related to a quadratic polynomial. These results will be useful in section \ref{ADE}. Moreover their proof, based on an induction process and on the additivity of the virtual Poincar\'e polynomial, will be the prototype for the computations of virtual Poincar\'e polynomial of spaces of arcs in this paper.

Let $p,q \in \mathbb N \cup 0$.
Let's denote by $Q_{p,q}$, or simply $Q$ depending on the context, the quadratic polynomial:
$$Q_{p,q}(y)=\sum_{i=1}^{p}y_i^2-\sum_{j=1}^{q}y_{j+p}^2.$$ 
First we recall the virtual Poincar\'e polynomial of the algebraic set $Y_{p,q}$ and $Y_{p,q}^{\epsilon}$ defined by
$Y_{p,q}=\{Q_{p,q}(y)= 0\}$ and $Y_{p,q}^{\epsilon}=\{Q_{p,q}(y)=\epsilon\}$
for $\epsilon \in \{1,-1\}$.

\begin{prop}(\cite{fichou2})\label{rappel} 
\begin{enumerate}
\item $\beta(Y_{p,q})=u^{p+q-1}-u^{\max \{p,q\}-1}+u^{\min \{p,q\}}.$
\item If $p \leq q$, then $\beta(Y_{p,q}^1)=u^{q-1}(u^p-1)$.
\item If $p > q$, then $\beta(Y_{p,q}^1)=u^{q}(u^{p-1}+1)$.
\end{enumerate}
\end{prop}

We will also use the notations $Y_{p,q}^*=Y_{p,q}\setminus \{0\}$ and $Y_{p,q}^c=\mathbb R^{p+q}\setminus Y_{p,q}$.

Next proposition collects the virtual Poincar\'e polynomial of the spaces of arcs associated with $Q$. We will use it later to derive the virtual Poincar\'e polynomial of spaces of arcs associated to the $ADE$-singularities.

\begin{prop}\label{quadra} Let $\epsilon \in \{1,-1\}$ and $l>2$.
\begin{enumerate}
\item If $p=0$ then $\beta(A_{l}^1(Q))=0$ whereas if $q=0$ then $\beta(A_{l}^{-1}(Q))=0$.
\item If $p=q=1$, then for $l=2n+1$
$$\beta(A_{2n+1}^{\epsilon}(Q))=nu^{2n+1}\beta(Y_{1,1}^*)$$
whereas for $l=2n$
$$\beta(A_{2n}^{\epsilon}(Q))=(n-1)u^{2n}\beta(Y_{1,1}^*)+u^{2n}\beta(Y_{1,1}^{\epsilon}).$$
\item In the other cases, if $l$ is odd, say $l= 2n+1$, then
$$\beta(A_{2n+1}^{\epsilon}(Q))=u^{(n+1)(p+q)-1}\beta(Y_{p,q}^*)\frac{u^{n(p+q-2)}-1}{u^{(p+q-2)}-1},$$
whereas if $l$ is even, say $l= 2n$, then
$$\beta(A_{2n}^{\epsilon}(Q))=u^{(n+1)(p+q)-2}\beta(Y_{p,q}^*)\frac{u^{(n-1)(p+q-2)}-1}{u^{(p+q-2)}-1}+u^{n(p+q)+2n}\beta(Y_{p,q}^{\epsilon}).$$
\end{enumerate}
\end{prop}

Let us introduce some notation. We consider polynomial arcs $\gamma (t)$ of the form
$$\gamma (t)=(a_1t+\cdots+a_lt^l,c_1^1t+\cdots+c_l^1t^l,\ldots,c_1^{p+q}t+\cdots+c_l^{p+q}t^l)$$
and denote by $c_i$ the vector $(c_i^1,\ldots,c_i^{p+q})$. We denote by $\Phi_{p,q}$, or simply $\Phi$ if the context is clear, the function defined on $\mathbb R^{p+q}\times \mathbb R^{p+q}$ by
$$\Phi_{p,q}(x,y)=2\sum_{j=1}^px_jy_j-2\sum_{j=1}^qx_{p+j}y_{p+j}.$$

\begin{proof}
We focus on the general case, so that we assume $pq \neq 0$ and $(p,q) \neq (1,1)$. Let us assume $\epsilon =1$. We begin with the case $l=2n+1$.

The composition with $Q$ of an arc $\gamma(t)$ produces the formal series
$$Q(c_1)t^2+\Phi(c_1,c_2)t^3+\cdots+(Q(c_n)+\sum_{s=1}^{n-1} \Phi(c_s,c_{2n-s}))t^{2n}+(\sum_{s=1}^{n}\Phi(c_s,c_{2n+1-s}))t^{2n+1}+\cdots$$

 therefore the set $A_{2n+1}^1(Q)$ is the algebraic set defined by the system:
\begin{displaymath}
\left \{ \begin{array}{lllll}
Q(c_1) =0 \\
\Phi(c_1,c_2) =0 \\
\ldots\\
Q(c_n)+\sum_{s=1}^{n-1} \Phi(c_s,c_{2n-s}) =0 \\
\sum_{s=1}^{n}\Phi(c_s,c_{2n+1-s}) =1.
\end{array} \right.
\end{displaymath}

We compute its virtual Poincar\'e polynomial using the additivity property via a suitable decomposition of $A_{2n+1}^1(Q)$. Namely, we decompose $A_{2n+1}^1(Q)$ following the vanishing of $c_1^1,c_1^2,\ldots,c_1^p$ in order to make appear an induction process. 

First, if $c_1^1 \neq 0$, we must choose $c_1^i$, for $i=1,\ldots,p+q$, in $Y_{p,q} \setminus Y_{p-1,q}$. Now, we may choose freely all other variables in the system, except $c_2^1,c_3^1,\ldots,c_{2n}^1$ for which we impose the value in order the equalities in the system to be satisfied. This part is linearly isomorphic to $Y_{p,q} \setminus Y_{p-1,q} \times \mathbb R^{2n(p-1)+2nq+1}$, therefore its contribution to the virtual Poincar\'e polynomial of $A_{2n+1}^1(Q)$ is $u^{2n(p+q-1)+1}\beta(Y_{p,q} \setminus Y_{p-1,q})$.

When $c_1^1 =0$, assume first that $c_1^2 \neq 0$. We now choose the $c_1^i$ for $i=2,\ldots,p+q$ in $Y_{p-1,q} \setminus Y_{p-2,q}$, and similarly we may choose freely all other variables except $c_2^2,c_3^2,\ldots,c_{2n}^2$. This contribution is equal to $u^{2n(p+q-1)+1}\beta(Y_{p-1,q} \setminus Y_{p-2,q})$. 

We repeat the argument until $c_1^p \neq 0$. It produces a contribution equal to $u^{2n(p+q-1)+1}\beta(Y_{1,q}^*)$. Finally if $c_1^p = 0$, then $c_1^{j+p}=0$ for $j=1,\ldots,q$ so that the first step of the computation is completed with a contribution equals to $u^{2n(p+q)+1}\beta(Y_{p,q}^*)$ by additivity of the virtual Poincar\'e polynomial. 

The new system is similar to that defining $A_{2n-1}^1(Q)$, with new variables $\tilde c_s^i=c_{s+1}^i$ for $ i=1,\ldots,p+q$ and $s=1,\ldots,2n-1$. Note however that the remaining $p+q$ variables $c_{2n+1}^i$ for $ i=1,\ldots,p+q$ remain free. Therefore
$$\beta(A_{2n+1}^1(Q))=u^{2n(p+q-1)+1}\beta(Y_{p,q}^*)+u^{p+q}\beta(A_{2n-1}^1(Q)).$$

An easy induction implies
$$\beta(A_{2n+1}^1(Q))=\beta(Y_{p,q}^*)\sum_{s=1}^{n}u^{(2n+1-s)(p+q-1)+s}=u^{(n+1)(p+q)-1}\beta(Y_{p,q}^*)\frac{u^{n(p+q-2)}-1}{u^{p+q-2}-1}$$
since at the last step, namely when $c_n^i=0$ for $i=1,\ldots,p+q$, the remaining system does not admit a solution.

In case $l=2n$ is even, the same method applies with the difference that at the last step the remaining equation, namely $Q(c_n)=1$, still admits solutions.
As a consequence
$$\beta(A_{2n}^1(Q))=\beta(Y_{p,q}^*)\sum_{s=1}^{n-1}u^{(2n-s)(p+q-1)+s}+u^{n(p+q)}\beta(Y_{p,q}^1).$$
\end{proof}
%%%%%%%%%%%%%%%%%%%%%%%%%%%%%%%%%%%%%%%%%%%%%%%%%%%%%%%%%%%%%%%%%%%%%%%%%%%%%%%%%%%%%%%%%%%%%%%%%%%%%%%%%%%%%%
\section{Blow-Nash types of ADE-singularities}\label{ADE}

We know that the corank and index of a Nash function germ are invariant under blow-Nash equivalence (Corollary \ref{invcorank}). In order to go further in the classification of singularities, the next step is to deal with simple singularities. Considering real analytic function germs, their simple singularities have been classified \cite{Arnold}: a real analytic function germ with a simple singularity is analytically equivalent to a polynomial germ belonging to one of the family:
\begin{displaymath}
\begin{array}{lllll}
A_k: \pm x^{k+1}+Q_{p,q}(y) \textrm{~~for ~~} k\geq 2,\\
D_k: x_1(\pm x_2^2 \pm x_1^{k-2})+Q_{p,q}(y) \textrm{~~for ~~} k\geq 4,\\
E_6: x_1^3 \pm x_2^4+Q_{p,q}(y),\\
E_7: x_1^3+x_1x_2^3+Q_{p,q}(y),\\
E_8: x_1^3+x_2^5+Q_{p,q}(y), \\
\end{array}
\end{displaymath}
for some $p,q \in \mathbb N\cup 0$.
This classification holds for Nash function germs. Indeed, let $f$ and $g$ be analytically equivalent Nash function germs. Then there exists an analytic isomorphism $h$ such that $f=g\circ h$. By Nash Approximation Theorem \cite{Shiota}, there exists also a Nash isomorphism $\widetilde h$ such that $f=g\circ \widetilde h$. Therefore $f$ and $g$ are Nash equivalent.

In this section, we state the classification of simple Nash function germs with respect to blow-Nash equivalence.
%%%%%%%%%%%%%%%%%%%%%%%%%%%%%%%%%%%%%%%%%%%%%%%%%%%%%%%%%%%%%%%%%%%%%%%%%%%%%%%%%%%%%%%%%%%%%%%%%%%%%%%%%%%%%%
\subsection{Blow-Nash types of $A_k$-singularities}

The $A_k$-singularities are the singularities given by the germs $f_k^{\pm}:(\mathbb R^{p+q+1},0) \rightarrow (\mathbb R,0)$ defined by
$$f_k^{\pm}(x,y)=\pm x^{k+1}+Q_{p,q}(y)$$
for $k\geq 2$ and some $p,q \in \mathbb N \cup 0$. 
Note that for $k$ even $f_k^+(x,y)=f_k^-(-x,y)$ so that $f_k^+$ and $f_k^-$ are linearly equivalent, thus blow-Nash equivalent. In that case we will simply write $f_k$.

We perform now some computations of virtual Poincar\'e polynomials that will be useful to distinguish the blow-Nash type of $A_k$-singularities in Proposition \ref{propAk}.
\begin{lemma}\label{lem3}\begin{flushleft}\end{flushleft}
\begin{enumerate}
\item \label{lem3.1} $\beta(\{f_{2n-1}^+=0\})=\beta(Y_{p+1,q})$ and $\beta(\{f_{2n-1}^-=0\})=\beta(Y_{p,q+1})$.
\item \label{lem3.2} If $p\leq q$, then $$\beta(\{f_{2n-1}^{+}=1\}) \neq \beta(\{f_{2n-1}^{-}=1\})=u\beta(Y_{p,q}^1),$$ whereas if $p \geq q$ then 
$$u\beta(Y_{p,q}^{-1})=\beta(\{f_{2n-1}^{+}=-1\}) \neq \beta(\{f_{2n-1}^{-}=-1\}).$$
\end{enumerate}
\end{lemma}

\begin{proof} Let us blow-up the origin in the algebraic variety defined by $f_{2n-1}^{\pm}(x,y)=0$. In the chart given by $x=u, y_i=v_iu$ with $i=1,\ldots,p+q$, the blowing-up variety is defined by
$$u^{2}f_{2n-3}^{\pm}(u,v)=0$$
where $v=(v_1,\ldots,v_{p+q})$.
The blowing-up is an isomorphism out from the origin between the strict transform $f_{2n-3}^{\pm}=0$ and $f_{2n-1}^{\pm}=0$, therefore
$$\beta(\{f_{2n-1}^{\pm}=0\}\setminus \{0\})=\beta(\{f_{2n-3}^{\pm}=0\}\setminus \{0\}).$$
As a consequence 
$$\beta(\{f_{2n-1}^{\pm}=0\})=\beta(\{f_{2n-3}^{\pm}=0\})$$
by additivity of the virtual Poincar\'e polynomial. Repeating $n-2$ times more the blowing-up at the origin leads to the first result.

For the second result, let us assume $p\leq q$. Applying the change of variables $u_i=y_i+y_{i+p}$, $v_i=y_i-y_{i+p}$ for $i=1,\ldots,p$, the equation $ f_{2n-1}^{\pm}=1 $ becomes
$$\pm x^{2n}+\sum_{i=1}^p u_iv_i-\sum_{j=p+1}^q y_{j}^2=1.$$
In order to compute the virtual Poincar\'e polynomial of this variety, we decompose it with respect to the vanishing of the variables $u_i$ for $i=1,\ldots,p$. If $u_1\neq 0$, we impose the value of $v_1$ so that the equality is satisfied. This contribution equals $(u-1)u^{p+q-1}$. If $u_1=0$, we repeat the argument with $u_2\neq 0$, for a contribution of $(u-1)u^{p+q-2}$. Repeating this argument until the case $u_1=\ldots=u_p=0$ implies that $\beta(\{f_{2n-1}^{\pm}=1\})$ is equal to
$$(u-1)\sum_{i=1}^p u^{p+q-i}+u^p \beta(\{\pm x^{2n}-\sum_{j=p+1}^q y_{j}^2=1\})=u\beta(Y_{p,q}^1)+u^p \beta(\{\pm x^{2n}-\sum_{j=p+1}^q y_{j}^2=1\}).$$
To conclude, remark that the non-empty variety $\{x^{2n}-\sum_{j=p+1}^q y_{j+p}^2=1\}$ has a nonzero virtual Poincar\'e polynomial contrary to $\{-x^{2n}-\sum_{j=p+1}^q y_{j+p}^2=1\}$.
\end{proof}
%%%%%%%%%%%%%%%%%%%%%%%%%%%%%%%%%%%%%%%%%%%%%%%%%%%%%%%%%%%%%%%%%%%%%%%%%%%%%%%%%%%%%%%%%%%%%%%%%%%%%%%%%%%%%%%
\subsubsection{Spaces of arcs for $A_k$-singularities}
We compute the virtual Poincar\'e polynomial of some spaces of arcs $A_l^{\pm 1}(f_k^{\pm})$ associated with the function germs $f_k^{\pm}$. Namely, we consider those polynomial arcs $\gamma (t)$ defined by
$$f_k^{\pm}(\gamma(t))=\pm t^l+\cdots$$
where the dots stand for higher order terms.
In case $k$ is sufficiently large with respect to $l$, only the quadratic part of $f_k^{\pm}$ plays a role in the space of arcs $A_l^{\pm 1}(f_k^{\pm})$. Therefore the next lemma is a direct consequence of Proposition \ref{quadra}.

\begin{lemma}\label{lem1} Let $k \geq l>2$ and $\epsilon \in \{1,-1\}$. Then $\beta(A_l^{\epsilon}(f_k^{\pm}))=u^l\beta(A_l^{\epsilon}(Q_{p,q})).$
\end{lemma}

Second, we deal with the first space of arcs associated with $f_k^{\pm}$ that takes into account the $x$-variable, namely $A_{k+1}^{\epsilon}(f_k^{\pm})$.

\begin{lemma}\label{lem2} Take ${\epsilon} \in \{1,-1\}$. 
%\begin{flushleft}\end{flushleft} 
\begin{enumerate}\item \label{lem2.1} $\beta(A_{2n+1}^{\epsilon}(f_{2n}))=u^{(n+1)(p+q)+2n}+u^{2n+1}\beta(A_{2n+1}^{\epsilon}(Q)).$
\item \label{lem2.2} $\beta(A_{2n}^{\epsilon}(f_{2n-1}^{\pm}))=u^{2n}\beta(A_{2n}^1(Q))+u^{n(p+q)+2n-1}(\beta(\{f_{2n-1}^{\pm}={\epsilon}\})-u\beta(Y_{p,q}^{\epsilon})).$

\end{enumerate}
\end{lemma}

\begin{proof}
\begin{flushleft}\end{flushleft}
\begin{enumerate}
\item The system of equations is the same as that of Lemma \ref{quadra} up to the last equation that has an additional additive term ``$+a_1^{2n+1}$'', namely
$$a_1^{2n+1}+\sum_{s=1}^{n}\Phi(c_s,c_{2n+1-s})=\epsilon.$$
The same method applies, with the only difference that at the last step, namely when $c_1^i=\ldots=c_n^i=0$ for $i=1,\ldots,p+q$, the remaining equation $a_1^{2n+1}=\epsilon$ still admits a solution. Therefore the value of $A_{2n+1}^{\epsilon}(f_{2n})$ under the virtual Poincar\'e polynomial is the sum of that of $A_{2n+1}^{\epsilon}(Q)$ plus the last contribution that equals $u^{(n+1)(p+q)+2n}$.
\item Similarly, the last equation of the system changes with respect to that of Lemma \ref{quadra}, so that the system of equations defining $A_{2n}^{\epsilon}(f_{2n-1})$ is
\begin{displaymath}
\left\{ \begin{array}{lllll}
Q(c_1)=0\\
\Phi(c_1,c_2)=0\\
\ldots\\
\sum_{s=1}^{n-1}\Phi(c_s,c_{2n-1-s})=0\\
a_1^{2n}+Q(c_n)+\sum_{s=1}^{n-1}\Phi(c_s,c_{2n-s})=\epsilon.\\
\end{array} \right.
\end{displaymath}
Once more the first part of the computation is the same as that of Lemma \ref{quadra}, whereas at the last step the variable $a_1$ is no longer free. As a consequence the term $u^{n(p+q)+2n}\beta(Y_{p,q}^{\epsilon})$ is replaced by $u^{n(p+q)+2n-1}\beta(\{f^{\pm}_{2n-1}=\epsilon\})$.
\end{enumerate}
\end{proof}

%%%%%%%%%%%%%%%%%%%%%%%%%%%%%%%%%%%%%%%%%%%%%%%%%%%%%%%%%%%%%%%%%%%%%%%%%%%

%%%%%%%%%%%%%%%%%%%%%%%%%%%%%%%%%%%%%%%%%%%%%%%%%%%%%%%%%%%%%%%%%%%%%%%%%%%%%%
\subsubsection{Classification of $A_k$-singularities}
Recall that in case $k$ is even, the function germs $f_k^+$ and $f_k^-$ are blow-Nash equivalent. We distinguish all other blow-Nash types by using the invariance under the blow-Nash equivalence of the zeta functions recalled in Theorem \ref{inv}. Recall that equivalent function germs share the same quadratic part (Corollary \ref{invcorank}).
\begin{prop}\label{propAk} If $f_k^{\epsilon_k}$ is blow-Nash equivalent to $f_l^{\epsilon_l}$, with $\epsilon_k,\epsilon_l \in \{+,-\}$, then $k=l$. Moreover if $k=l$ is odd, then $\epsilon_k=\epsilon_l$.
\end{prop}

\begin{cor} Two germs with $A_k$-singularities are blow-Nash equivalent if and only if they are analytically equivalent.
\end{cor}

\begin{proof}[Proof of Proposition \ref{propAk}] If $k \neq l$ assume for instance $k<l$. For even $k$, then $\beta(A_{k+1}^1(f_k))$ is different from $\beta(A_{k+1}^1(f_l^{\pm}))$ thanks to Lemma \ref{lem2}.\ref{lem2.1}. For $k$ odd, Lemma \ref{lem2}.\ref{lem2.2} combined with Lemma \ref{lem3}.\ref{lem3.2} assert that either 
$$\beta(A_{k+1}^1(f_k^{\pm})) \neq \beta(A_{k+1}^1(f_l^{\pm}))$$
or 
$$\beta(A_{k+1}^{-1}(f_k^{\pm})) \neq \beta(A_{k+1}^{-1}(f_l^{\pm})).$$
Therefore the zeta functions of $f_k^{\pm}$ and $f_l^{\pm}$ are different, and thus $f_k^{\pm}$ and $f_l^{\pm}$ can not be blow-Nash equivalent.

If $k=l$ is odd, Lemma \ref{lem2}.\ref{lem2.2} combined with Lemma \ref{lem3}.\ref{lem3.2} assert once more that either 
$$\beta(A_{k+1}^1(f_k^{+})) \neq \beta(A_{k+1}^1(f_k^{-})),$$
in the case $p \leq q$, or 
$$\beta(A_{k+1}^{-1}(f_k^{+})) \neq \beta(A_{k+1}^{-1}(f_l^{-}))$$
if $p> q$, so that $f_k^{+}$ and $f_k^{-}$ can not be blow-Nash equivalent.
\end{proof}
%%%%%%%%%%%%%%%%%%%%%%%%%%%%%%%%%%%%%%%%%%%%%%%%%%%%%%%%%%%%%%%%%%%%%%%%%%%%%%%%%%%%%%%%%%%%%%
%%%%%%%%%%%%%%%%%%%%%%%%%%%%%%%%%%%%%%%%%%%%%%%%%%%%%%%%%%%%%%%%%%%%%%%%%%%%%%%%%%%%%%%%%%%%%%
%%%%%%%%%%%%%%%%%%%%%%%%%%%%%%%%%%%%%%%%%%%%%%%%%%%%%%%%%%%%%%%%%%%%%%%%%%%%%%%%%%%%%%%%%%%%%%
\subsection{Blow-Nash types of $D_k$-singularities}\label{sectionD}

The $D_k$-singularities are the function germs $g_k^{\epsilon_1,\epsilon_2}:(\mathbb R^{p+q+2},0) \rightarrow (\mathbb R,0)$ defined by
$$g_k^{\epsilon_1,\epsilon_2}(x_1,x_2,y)=x_1(\epsilon_1 x_2^2+\epsilon_2 x_1^{k-2})+Q_{p,q}(y)$$
for $k\geq 4$, $\epsilon_1,\epsilon_2 \in \{+,-\}$ and some $p,q \in \mathbb N \cup 0$.

Note that for $k$ odd, the equality $g_k^{+,\epsilon_2}(-x_1,x_2,y)=g_k^{-,\epsilon_2}(x_1,x_2,y)$ holds so that $g_k^{+,\epsilon_2}$ and $g_k^{-,\epsilon_2}$ are linearly equivalent, and therefore blow-Nash equivalent. Similarly, for $k$ even $g_k^{\epsilon_1,\epsilon_2}(-x_1,x_2,y)=g_k^{-\epsilon_1,-\epsilon_2}(x_1,x_2,y)$ so that $g_k^{\epsilon_1,\epsilon_2}$ and $g_k^{-\epsilon_1,-\epsilon_2}$ are blow-Nash equivalent.

The next lemma will be useful in order to determine the blow-Nash types of $D_k$-singularities in Proposition \ref{propDk}.

\begin{lemma}\label{lem6} \begin{flushleft} \end{flushleft}
\begin{enumerate}
\item If $k$ is odd, then $\beta(\{x_1x_2^2+x_1^{k-1}=1\})=2u$.
\item\label{lem6.2} If $k$ is even, then $\beta(\{x_1x_2^2+x_1^{k-1}=1\})=u$ and  $\beta(\{x_1x_2^2-x_1^{k-1}=1\})=2u$.
\end{enumerate}
\end{lemma}

\begin{proof} We compute the virtual Poincar\'e polynomial of these plane curves by compactifying the curve, and then resolving its singularities. Assume $k$ is odd. Compactifying the nonsingular curve $\{x_1x_2^2+x_1^{k-1}=1\}$ in the projective plane gives a singular curve $C$ with one point $p$, singular, at infinity. By invariance under algebraic isomorphisms of the virtual Poincar\'e polynomial, we obtain
$$\beta(\{x_1x_2^2+x_1^{k-1}=1\})=\beta(C\setminus \{p\}).$$

Let $\widetilde C$ be a resolution of the singularities of $C$.
Then $\widetilde C$ has two connected components homeomorphic to a circle, each of them containing a preimage, say $p_1$ and $p_2$, of $p$. Therefore $$\beta(C\setminus \{p\})=\beta(\widetilde C \setminus \{p_1,p_2\})=2(u+1)-2$$
by additivity of the virtual Poincar\'e polynomial. Finally
$$\beta(\{x_1x_2^2+x_1^{k-1}=1\})=\beta(C\setminus \{p\})=2u.$$

The results in the case $k$ even come from the fact that the resolution of the singularities of the nonsingular plane curve $\{x_1x_2^2+x_1^{k-1}=1\}$ has one component, whereas that of $\{x_1x_2^2-x_1^{k-1}=1\}$ has two components.
\end{proof}
%%%%%%%%%%%%%%%%%%%%%%%%%%%%%%%%%%%%%%%%%%%%%%%%%%%%%%%%%%%%%%%%%%%%%%%%%%%%%%%%%%%%%%%%%%%%%%%
%%%%%%%%%%%%%%%%%%%%%%%%%%%%%%%%%%%%%%%%%%%%%%%%%%%%%%%%%%%%%%%%%%%%%%%%%%%%%%%%%%%%%%%%%%%%%%%
\subsubsection{Spaces of arcs for $D_k$-singularities}
Similarly to the $A_k$-singularity case, we compute separately the virtual Poincar\'e polynomial of the spaces of arcs $A_{k-1}^{\pm 1}(g_k)$ and  $A_{l}^{\pm 1}(g_k)$ where $k > l+1$.

In order to fix some notation, we consider those polynomial arcs $\gamma (t)$ of the form
$$(a_1t+\cdots+a_lt^l,b_1t+\cdots+b_lt^l,c_1^1t+\cdots+c_l^1t^l,\ldots,c_1^{p+q}t+\cdots+c_l^{p+q}t^l)$$
defined by
$$g_k^{\epsilon_1,\epsilon_2}(\gamma(t))=\pm t^l+\cdots$$

As an intermediate step, we are interested in the spaces of arcs for the function $G(x_1,x_2,y)=x_1 x_2^2+Q_{p,q}(y)$.

\begin{lemma}\label{lem4} Let $l \geq 3$ be an integer and take $\epsilon \in \{1,-1\}$.
%\begin{flushleft}\end{flushleft}
\begin{enumerate}
\item If $(p,q)=(0,0)$ then \begin{displaymath} \beta(A_{l}^{\epsilon}(G))=
\left\{ \begin{array}{ll}
u^{2n+2}(u^{n}-1) \textrm{~~if~~} l=2n+1\\
u^{2n+1}(u^{n-1}-1) \textrm{~~if~~} l=2n.\\
\end{array} \right.
\end{displaymath}

\item When $l$ is odd, say $l=2n+1$, if $(p,q)=(1,0)$ or $(p,q)=(0,1)$ then $$\beta(A_{2n+1}^{\epsilon}(G))=n(u-1)u^{4n+2}.$$
\item When $l$ is even, say $l=2n$, and $p+q=1$ then
\begin{displaymath}
\beta(A_{2n}^{\epsilon}(G))=
\left\{ \begin{array}{ll}
(n-1)(u-1)u^{4n}+2u^{4n+1} \textrm{~~if~~} (\epsilon,p) =(1,1)\textrm{~~or~~}(\epsilon,p) =(-1,0),\\
(n-1)(u-1)u^{4n} \textrm{~~if~~} (\epsilon,p) =(-1,1)\textrm{~~or~~}(\epsilon,p) =(1,0).\\
\end{array} \right.
\end{displaymath}

\item When $l$ is odd, say $l=2n+1$, then $\beta(A_{2n+1}^{\epsilon}(G))$ is equal to
$$u^{(n+1)(p+q)}(\frac{u^{n(p+q-1)}-1}{u^{p+q-1}-1}u^{5}\beta(Y_{p,q}^*)
+\frac{u^{(n-1)(p+q-1)}-1}{u^{p+q-1}-1}(u-1)u^{(p+q)+3}+u^{3n+1}(u-1)).$$
\item When $l$ is even, say $l=2n$, then
$$\beta(A_{2n}^{\epsilon}(G))= \frac{u^{(n-1)(p+q-1)}-1}{u^{p+q-1}-1}u^{(n+2)(p+q)+3n} (\beta(Y_{p,q}^*)+(u-1))
+u^{n(p+q)+3n+1}\beta(Y_{p,q}^{\epsilon}).$$
\end{enumerate}
\end{lemma}

\begin{proof} We concentrate on the general case $p+q >1$.

(4) The system of equations defining $A_{2n+1}^{\epsilon}(G)$ is:
\begin{displaymath}
\left\{ \begin{array}{llllll}
Q(c_1)=0\\
a_1b_1^2+\Phi(c_1,c_2)=0\\
2a_1b_1b_2+a_2b_1^2+Q(c_2)+\Phi(c_1,c_3)=0\\
\ldots\\
\sum_{s=1}^{n-1} b_s^2a_{2n-2s} +2\sum_{s=1}^{n-2}b_s\sum_{t=2}^{2n-2s} b_ta_{2n-t-s}+Q(c_n)+\sum_{s=1}^{n-1}\Phi(c_s,c_{2n-s})=0\\
\sum_{s=1}^{n} b_s^2a_{2n+1-2s} +2\sum_{s=1}^{n-1}b_s\sum_{t=2}^{2n+1-2s} b_ta_{2n+1-t-s}+\sum_{s=1}^{n}\Phi(c_s,c_{2n+1-s})=\epsilon.\\
\end{array} \right.
\end{displaymath}

We adapt the proof of Lemma \ref{quadra} to the case of $G$ in place of the quadratic polynomial $Q$.
 At the first step, we consider the disjunction of cases given by $c_1^1 \neq 0$, then by $c_1^1 =0$ and $c_1^2 \neq 0$, $\ldots$, finally by $c_1^1=\ldots=c_1^{p-1}=0$ and $c_1^p \neq 0$. Note that the $a_s$ and $b_s$ variables do not intervene. So, similarly to the proof of Lemma \ref{quadra}, the corresponding contribution is $u^{2n(p+q-1)+2(2n+1)}\beta(Y_{p,q}^*)$.

Now, when $c_1^1=\ldots=c_1^{p}=0$, then $c_1^{p+1},\ldots,c_1^{p+q}$ must vanish and the second equation imposes the condition $a_1b_1^2=0$. If $b_1 \neq 0$, then necessarily $a_1=0$. We may describe the solutions of the system by imposing the value of $a_2,\ldots,a_{2n-1}$. This part contributes to $\beta(A_{2n+1}^1(G))$ as $(u-1)u^{2n+2+2n(p+q)}$.

The remaining equations correspond to the system defining $A_{2n-1}^{\epsilon}(G)$ via the change of coordinates $\tilde a_s=a_s, \tilde b_s=b_{s+1}, \tilde c_s^i =c_{s+1}^i$ with $s=1,\ldots,2n-1$ and $i=1,\ldots,p+q$. Note that the variables $a_{2n},a_{2n+1},b_{2n+1},c_{2n+1}^i$, for $i=1,\ldots,p+q$ are free. Therefore
$$\beta(A_{2n+1}^{\epsilon}(G))=u^{2n(p+q-1)+2(2n+1)}\beta(Y_{p,q}^*)+(u-1)u^{2n+2+2n(p+q)}+u^{p+q+3}\beta(A_{2n-1}^{\epsilon}(G)).$$
We repeat this argument until obtaining $\beta(A_3^{\epsilon}(G))$. A similar computation implies
$$\beta(A_3^{\epsilon}(G))=u^{2(p+q)+5}\beta(Y_{p,q}^*)+u^{2(p+q)+4}\beta(\{x_1x_2^2=\epsilon\}).$$
As $\beta(\{x_1x_2^2=\epsilon\})$ is simply equal to $u-1$, it follows that $\beta(A_{2n+1}^{\epsilon}(G))$ equals
\begin{displaymath}
\begin{array}{ll} u^{(n+2)(p+q)+3n}\beta(Y_{p,q}^*)\sum_{s=0}^{n-1} u^{s(p+q-1)}
+(u-1)u^{(n+3)(p+q)+3n-1}\sum_{s=0}^{n-2} u^{s(p+q-1)}\\
+u^{(n+1)(p+q)+3n+1}(u-1).
\end{array}
\end{displaymath}

%\item 

(5) The same method applies to compute the virtual Poincar\'e polynomial of $A_{2n}^{\epsilon}(G)$. In this way we obtain the following relations between the virtual Poincar\'e polynomial of the spaces of arcs $A_{2n}^{\epsilon}(G)$ and $A_{2n-2}^{\epsilon}(G)$:
$$\beta(A_{2n}^{\epsilon}(G))=u^{(2n-1)(p+q-1)+4n}\beta(Y_{p,q}^*)+(u-1)u^{2n+1+(2n-1)(p+q)}+u^{p+q+3}\beta(A_{2n-2}^{\epsilon}(G)).$$
Repeating the argument $n-1$-times leads to the virtual Poincar\'e polynomial of $A_2^{\epsilon}(G)$. It equals  
$$\beta(A_2^{\epsilon}(G))=u^{4+p+q}\beta(Y_{p,q}^{\epsilon}).$$
Finally the formula computing $\beta(A_{2n}^{\epsilon}(G))$ is:
\begin{displaymath}
\begin{array}{ll} u^{(n+2)(p+q)+3n}\beta(Y_{p,q}^*)\sum_{s=0}^{n-2} u^{s(p+q-1)}
+(u-1)u^{(n+3)(p+q)+3n-1}\sum_{s=0}^{n-2} u^{s(p+q-1)}\\
+u^{n(p+q)+3n+1}\beta(Y_{p,q}^{\epsilon}).
\end{array}
\end{displaymath}
%\end{enumerate}
\end{proof}

In case $k>l+1$, the spaces of arcs of $g_k^{\epsilon_1,\epsilon_2}$ and $G$ are similar since the monomial $\pm x_1^{k-1}$ has a higher degree.

\begin{cor} Take $l\geq 3$. If $k>l+1$ then $\beta(A_l^{\epsilon}(g_k^{\epsilon_1,\epsilon_2}))=\beta(A_l^{\epsilon}(G))$ for $\epsilon \in \{1,-1\}$ and $\epsilon_1,\epsilon_2 \in \{+,-\}$.
\end{cor}

We compute now the virtual Poincar\'e polynomial of $A_{k+1}^{\pm 1}(g_k)$.

\begin{lemma}\label{lem5} Take $\epsilon \in \{1,-1\}$, $\epsilon_1,\epsilon_2 \in \{+,-\}$ and $n\geq 1$. 
%\begin{flushleft} \end{flushleft}
\begin{enumerate}
\item\label{lem5bis} If $(p,q)=(0,0)$ and $k\geq 4$ then
\begin{displaymath}
\beta(A_{k-1}^{\epsilon}(g_{k}^{\epsilon_1,\epsilon_2}))=\beta(A_{k-1}^{\epsilon}(G))+ \left \{
\begin{array}{ll} 
u^{3n-1} \textrm{~~if~~} k=2n, \\
u^{3n}\beta(\{g_k^{\epsilon_1,\epsilon_2}(x_1,0,0)=\epsilon\}) \textrm{~~if~~} k=2n+1.\\
\end{array} \right.
\end{displaymath}

\item\label{lem5.1} $\beta(A_{2n+1}^{\epsilon}(g_{2n+2}^{\epsilon_1,\epsilon_2}))$ is equal to 
$$\beta(A_{2n+1}^{\epsilon}(G))+u^{(n+1)(p+q)+3n+1}(\beta(\{g_{2n+2}^{\epsilon_1,\epsilon_2}(x_1,x_2,0)=\epsilon\})-(u-1)).$$

\item\label{lem5.2} $\beta(A_{2n}^{\epsilon}(g_{2n+1}^{\epsilon_1,\epsilon_2}))$ is equal to
$$\beta(A_{2n}^{\epsilon}(G))+u^{n(p+q)+3n}(\beta(\{g_{2n+1}^{\epsilon_1,\epsilon_2}(x_1,0,y)=\epsilon\})-u\beta(Y_{p,q}^{\epsilon})).$$
\end{enumerate}
\end{lemma}

\begin{proof} Let us focus on point $(2)$, the proof of points $(1)$ and $(3)$ being similar. Computations are analog to that of Lemma \ref{lem4}, except concerning the last equation of the system defining $A_{2n+1}^{\epsilon}(g_{2n+2}^{\epsilon_1,\epsilon_2})$ in which an additional additive term $a_1^{2n+1}$ appears. Therefore, following the induction process of the proof of Lemma \ref{lem4}, the preceding contribution of $\beta(A_3^{\epsilon}(G))$ coming from the last equation, namely 
$$\beta(A_3^{\epsilon}(G))=u^{2(p+q)+5}\beta(Y_{p,q}^*)+u^{2(p+q)+4}\beta(\{x_1x_2^2=\epsilon\}),$$
 is replaced by
$$\beta(A_3^{\epsilon}(G))=u^{2(p+q)+5}\beta(Y_{p,q}^*)+u^{2(p+q)+4}\beta(\{g_{2n+2}^{\epsilon_1,\epsilon_2}(x_1,x_2,0)=\epsilon\}).$$
\end{proof}

\begin{rmk}\label{rmklem5} Note that the computations of the virtual Poincar\'e polynomial related to the spaces of arcs for the naive zeta function are similar. For instance $\beta(A_{2n}(g_{2n+1}^{\epsilon_1,\epsilon_2}))$ is equal to
$$\beta(A_{2n}(G))+u^{n(p+q)+3n}(\beta(\{g_{2n+1}^{\epsilon_1,\epsilon_2}(x_1,0,y)\neq 0\})-u\beta(Y_{p,q}^{c})).$$
\end{rmk}
%%%%%%%%%%%%%%%%%%%%%%%%%%%%%%%%%%%%%%%%%%%%%%%%%%%%%%%%%%%%%%%%%%%%%%%%%%%%%%%%%%%%%%%%%%%%%

%%%%%%%%%%%%%%%%%%%%%%%%%%%%%%%%%%%%%%%%%%%%%%%%%%%%%%%%%%%%%%%%%%%%%%%%%%%%%%%%%%%%%%%%%%%%%%%
\subsubsection{Classification of $D_k$-singularities}
Take $\epsilon_1,\epsilon_2 \in \{+,-\}$.
Recall that in the case $k$ odd, the germs $g_k^{+,\epsilon_2}$ and $g_k^{-,\epsilon_2}$ are blow-Nash equivalent whereas in the case $k$ even $g_k^{\epsilon_1,\epsilon_2}$ and $g_k^{-\epsilon_1,-\epsilon_2}$ are blow-Nash equivalent.

\begin{prop}\label{propDk} Take $\epsilon_1,\epsilon_2,\sigma_1,\sigma_2  \in \{+,-\}$. If $g_k^{\epsilon_1,\epsilon_2}$ is blow-Nash equivalent to $g_l^{\sigma_1,\sigma_2}$, then $k=l$. If moreover $k=l$ is odd then $\epsilon_2=\sigma_2$, whereas if $k$ is even then $\epsilon_1\epsilon_2=\sigma_1\sigma_2$.
\end{prop}

\begin{cor} Two germs with $D_k$-singularities are blow-Nash equivalent if and only if they are analytically equivalent.
\end{cor}

\begin{proof}[Proof of Proposition \ref{propDk}] Assume $k<l$. We prove that the coefficients of $T^{k-1}$ in at least one of the zeta functions of $g_k^{\epsilon_1,\epsilon_2}$ and $g_l^{\sigma_1,\sigma_2}$ are different.

If $(p,q)=(0,0)$ the result comes directly from Lemma \ref{lem5}.\ref{lem5bis}. Otherwise, if $k$ is even it suffices to prove that
$$\beta(\{g_{k}^{\epsilon_1,\epsilon_2}(x_1,x_2,0)=\epsilon\})-(u-1)$$
can not be the zero polynomial thanks to Lemma \ref{lem5}.\ref{lem5.1}. This is guaranteed by Lemma \ref{lem6}.\ref{lem6.2}.

If $k$ is odd, it suffices to prove by Remark \ref{rmklem5} that
$$\beta(\{g_{k}^{\epsilon_1,\epsilon_2}(x_1,0,y)\neq 0\})-u\beta(Y_{p,q}^{c})$$
can not be the zero polynomial. If $\epsilon_2 =+$ note that
$$\beta(\{g_{k}^{\epsilon_1,+}(x_1,0,y)\neq 0\})=\beta(Y_{p+1,q}^{c})$$
by Lemma \ref{lem3}.\ref{lem3.1}. But $\beta(Y_{p+1,q})$ can not be equal to $u\beta(Y_{p,q})$ thanks to Lemma \ref{rappel}. Actually for $p<q$ the former equals $u^{p+q}-u^{q-1}+u^{p+1}$ whereas the latter is $u^{p+q}-u^{q}+u^{p+1}$, and for $p=q$ the former equals $u^{2p}$ whereas the latter is $u^{2p}-u^{p}+u^{p+1}$. Therefore the $T^{k-1}$ coefficients of the naive zeta functions of $g_k^{\epsilon_1,+}$ and $g_l^{\sigma_1,\sigma_2}$ are different. Similarly if $\epsilon_2 =-$ then 
$$\beta(\{g_{k}^{\epsilon_1,-}(x_1,0,y)\neq 0\})=\beta(Y_{p,q+1}^{c})$$
which can not be equal to $u\beta(Y_{p,q}^c)$. Thus the $T^{k-1}$ coefficients of the naive zeta functions of $g_k^{\epsilon_1,-}$ and $g_l^{\sigma_1,\sigma_2}$ are still different.

%If $k$ is odd, it suffices to prove by Lemma \ref{lem5}.\ref{lem5.2} that
%$$\beta(\{g_{k}^{\epsilon_1,\epsilon_2}(x_1,0,y)=\epsilon\})-u\beta(Y_{p,q}^{\epsilon})$$
%can not be the zero polynomial. If $\epsilon_2 =+$ note that
%$$\beta(\{g_{k}^{\epsilon_1,+}(x_1,0,y)=\epsilon\})=\beta(Y_{p+1,q}^{\epsilon}).$$
% But $\beta(Y_{p+1,q}^{1})$ can not be equal to $u\beta(Y_{p,q}^{1})$ thanks to Lemma \ref{rappel}. Therefore the $T^{k-1}$ coefficients of $Z^{1}_{g_k^{\epsilon_1,\epsilon_2}}$ and $Z^{1}_{g_l^{\sigma_1,\sigma_2}}$ are different. Similarly if $\epsilon_2 =-$ the $T^{k-1}$ coefficients of $Z^{-1}_{g_k^{\epsilon_1,\epsilon_2}}$ and $Z^{-1}_{g_l^{\sigma_1,\sigma_2}}$ are different.

If $k=l$ is odd, the difference between the $T^{k-1}$-coefficients in the zeta functions with sign of $g_k^{\epsilon_1,\epsilon_2}$ and $g_k^{\sigma_1,\sigma_2}$ may appear in the expressions $\beta(\{\epsilon_2x_1^{k-1}+Q(y)=\epsilon\})$ and $\beta(\{\sigma_2x_1^{k-1}+Q(y)=\epsilon\})$ by Lemma \ref{lem5}.\ref{lem5.2}. But $k-1$ being even, these polynomials are equal if and only if $\epsilon_2=\sigma_2$ by Lemma \ref{lem3}.\ref{lem3.2}. The argument is similar for $k=l$ even, combining Lemma \ref{lem5}.\ref{lem5.1} and Lemma \ref{lem6}. 
\end{proof}
%%%%%%%%%%%%%%%%%%%%%%%%%%%%%%%%%%%%%%%%%%%%%%%%%%%%%%%%%%%%%%%%%%%%%%%%%%%%%%%%%%%%%%%%%%%%%
%%%%%%%%%%%%%%%%%%%%%%%%%%%%%%%%%%%%%%%%%%%%%%%%%%%%%%%%%%%%%%%%%%%%%%%%%%%%%%%%%%%%%%%%%%%%%
%%%%%%%%%%%%%%%%%%%%%%%%%%%%%%%%%%%%%%%%%%%%%%%%%%%%%%%%%%%%%%%%%%%%%%%%%%%%%%%%%%%%%%%%%%%%%
\subsection{Blow-Nash types of $E_6,E_7,E_8$-singularities}
%%%%%%%%%%%%%%%%%%%%%%%%%%%%%%%%%%%%%%%%%%%%%%%%%%%%%%%%%%%%%%%%%%%%%%%%%%%%%%%%%%%%%%%%%%%%%

The $E_6,E_7,E_8$-singularities are the function germs defined on $(\mathbb R^{p+q+2},0)$ by
\begin{displaymath}
\begin{array}{lll}
h_6^{\pm}(x_1,x_2,y)=x_1^3 \pm x_2^{4}+Q_{p,q}(y),\\
h_7(x_1,x_2,y)=x_1^3 + x_1 x_2^{3}+Q_{p,q}(y),\\
h_8(x_1,x_2,y)=x_1^3 + x_2^{5}+Q_{p,q}(y),\\
\end{array}
\end{displaymath}
for some $p,q \in \mathbb N \cup 0$.

%%%%%%%%%%%%%%%%%%%%%%%%%%%%%%%%%%%%%%%%%%%%%%%%%%%%%%%%%%%%%%%%%%%%%%%%%%%%%%%%%%%%%%%%%%%%%
\subsubsection{Spaces of arcs for $E_6,E_7,E_8$-singularities}
We compute the virtual Poincar\'e polynomial of some spaces of arcs related to $h_6^+,h_6^-,h_7,h_8$.

\begin{lemma}\label{lem7} Let $\epsilon \in \{1,-1\}$.%\begin{flushleft} \end{flushleft}
\begin{enumerate}
\item\label{lem7.1} For $h=h_6^+,h_6^-,h_7$ or $h_8$, then
$\beta(A_3^{\epsilon}(h))=u^{2(p+q)+7}\beta(Y_{p,q}^*)+u^{2(p+q)+5}.$
\item\label{lem7.2} \begin{displaymath} \beta(A_4(h_6^{\sigma}))=
\left\{ \begin{array}{ll}
(u-1)u^{3(p+q)+6}\beta(Y_{p,q}^*)+u^{2(p+q)+6}\beta(Y_{p+1,q}^c) \textrm{~~if ~~} \sigma=+,\\
(u-1)u^{3(p+q)+6}\beta(Y_{p,q}^*)+u^{2(p+q)+6}\beta(Y_{p,q+1}^c) \textrm{~~if ~~} \sigma=-,\\
\end{array} \right.
\end{displaymath}
and
$$\beta(A_4(h_7))=\beta(A_4(h_8))=(u-1)u^{3(p+q)+6}\beta(Y_{p,q}^*)+u^{2(p+q)+7}\beta(Y_{p,q}^c).$$
\item\label{lem7.suppl} $\beta(A_4^{\epsilon}(h_7))=\beta(A_4^{\epsilon}(h_8))=u^{3(p+q)+6}\beta(Y_{p,q}^*)+u^{2(p+q)+7}\beta(Y_{p,q}^{\epsilon}).$
\item\label{lem7.3} \begin{displaymath} \beta(A_5^{\epsilon}(h_i))=
\left\{ \begin{array}{ll}
\beta(Y_{p,q}^*)(u^{4(p+q)+7}+u^{3(p+q)+8})+(u-1)u^{3(p+q)+7} \textrm{~~if ~~} i=7,\\
\beta(Y_{p,q}^*)(u^{4(p+q)+7}+u^{3(p+q)+8})+u^{3(p+q)+8} \textrm{~~if ~~} i=8.\\
\end{array} \right.
\end{displaymath}
\end{enumerate}
\end{lemma}

%The proof is very similar to the previous ones. Nevertheless, note that for the computation of $\beta(A_4(h_6^{\pm}))$ we make use of lemma \ref{lem3}.\ref{lem3.1}.

\begin{proof} The proof of (\ref{lem7.1}),(\ref{lem7.suppl}) and (\ref{lem7.3}) run as that of Lemma \ref{lem2}. For (\ref{lem7.2}), let us deal with $h_6^+$ for example (the other cases are similar). The system of equations defining $A_4(h_6^{+})$ is
\begin{displaymath}
\left\{ \begin{array}{lll}
Q(c_1)=0\\
a_1^3+ \Phi(c_1,c_2)=0\\
3a_1^2a_2+b_1^4+Q(c_2)+\Phi(c_1,c_3)\neq 0.\\
\end{array} \right.
\end{displaymath}
We use the following decomposition of $A_4(h_6^{+})$ in order to compute its virtual Poincar\'e polynomial. Let's choose a nonzero $p+q$-uplet $(c_1^1,\ldots,c_1^{p+q}) \in Y_{p,q}^*$ satisfying the first equation. Say for example $c_1^1\neq 0$. To describe the solutions of the system, we may impose the choice of $c_2^1$ in order the second equation to be satisfied, and exclude one value for $c_3^1$ for the last equation. This part contributes to $\beta(A_4(h_6^{+}))$ as $(u-1)u^{3(p+q)+6}\beta(Y_{p,q}^*)$. Now if $c_1^1,\ldots,c_1^{p+q}$ vanish then $a_1$ must vanish also and the last equation becomes
$$b_1^4+Q(c_2) \neq 0.$$
The corresponding virtual Poincar\'e polynomial is $\beta(Y_{p+1,q}^c)$ by Lemma \ref{lem3}.\ref{lem3.1}, hence an additional contribution of $u^{2(p+q)+6}\beta(Y_{p+1,q}^c)$ that completes the proof.
\end{proof}
%%%%%%%%%%%%%%%%%%%%%%%%%%%%%%%%%%%%%%%%%%%%%%%%%%%%%%%%%%%%%%%%%%%%%%%%%%%%%%%%%%%%%%%%%%%%%
\subsubsection{Classification of $E_6,E_7,E_8$-singularities}
The classification of $E_6,E_7,E_8$-singularities under blow-Nash equivalence is obtained by combining Lemma \ref{lem7} and the invariance of the zeta functions of Theorem \ref{inv}.

\begin{prop}\label{propE} The function germs $h_6^+,h_6^-,h_7,h_8$ belong to different blow-Nash equivalence classes.
\end{prop}

\begin{cor} Two germs with $E_6^{\pm}, E_7$ or $E_8$-singularities are blow-Nash equivalent if and only if they are analytically equivalent.
\end{cor}

\begin{proof}[Proof of Proposition \ref{propE}] The germs $h_6^+$ and $h_6^-$ can not be blow-Nash equivalent by the first two equalities of Lemma \ref{lem7}.\ref{lem7.2}. It follows also from Lemma \ref{lem7}.\ref{lem7.2} that $h_6^+$ and $h_6^-$ belong to different blow-Nash equivalence classes from the classes of $h_7,h_8$. Finally $h_7$ and $h_8$ are not blow-Nash equivalent by Lemma \ref{lem7}.\ref{lem7.3}.
\end{proof}
%%%%%%%%%%%%%%%%%%%%%%%%%%%%%%%%%%%%%%%%%%%%%%%%%%%%%%%%%%%%%%%%%%%%%%%%%%%%%%%%%%%%%%%%%%%%%
%%%%%%%%%%%%%%%%%%%%%%%%%%%%%%%%%%%%%%%%%%%%%%%%%%%%%%%%%%%%%%%%%%%%%%%%%%%%%%%%%%%%%%%%%%%%%
%%%%%%%%%%%%%%%%%%%%%%%%%%%%%%%%%%%%%%%%%%%%%%%%%%%%%%%%%%%%%%%%%%%%%%%%%%%%%%%%%%%%%%%%%%%%%
\section{Classification of simple singularities}\label{sectionsimple}
We have proven so far that, inside each list of simple singularities, the blow-Nash classification coincides with the analytic one. Actually this result remains true without specifying the membership of any given list!

The following theorem states the classification under blow-Nash equivalence of Nash function germs equipped with a simple singularity.
 
\begin{thm}\label{thm1} Let $f,g:(\mathbb R^{d},0) \rightarrow (\mathbb R,0)$ be Nash function germs. Assume $f$ and $g$ are simple. Then $f$ and $g$ are blow-Nash equivalent if and only if $f$ and $g$ are analytically equivalent.
\end{thm}

The proof consists in comparing the different virtual Poincar\'e polynomials of the spaces of arcs associated with the $ADE$-singularities.

\begin{proof} A germ with a $A_k$-singularity can not be blow-Nash equivalent to a $D_l, E_6,E_7$ or $E_8$-singularity by Corollary \ref{invcorank} since the former has corank one whereas the latter has corank two. Therefore we are reduced to consider only corank two germs.

Assume $g=g_k$ has a $D_k$-singularity. The computations of Lemmas \ref{lem4} and \ref{lem7} prove that the virtual Poincar\'e polynomial of the space of arcs $A_3^1$ for $g_k$ is different from that of $h_6^{\pm},h_7,h_8$ if $k>4$. Therefore $g_k$ is not blow-Nash equivalent to $h_6^{\pm},h_7,h_8$ if $k>4$. It remains to prove that $g_4^{\pm}$ can not be blow-Nash equivalent to $h_6^{\pm},h_7$ or $h_8$. Note that $\beta(A_3^{\pm1}(g_4^{\pm}))$ is equal to the virtual Poincar\'e polynomial of the space of arcs $A_3^{\pm}$ for $h_6^{\pm},h_7,h_8$ by Lemma \ref{lem7}.

Actually:
\begin{displaymath}
\begin{array}{llll}
\beta(A_4^{\epsilon}(g_4^{+}))=u^{3(p+q)+6}\beta(Y_{p,q}^*)+(u-1)u^{3(p+q)+5}+u^{2(p+q)+6}\beta(Y_{p,q}^{\epsilon}),\\
\beta(A_4^{\epsilon}(g_4^{-}))=u^{3(p+q)+6}\beta(Y_{p,q}^*)+3(u-1)u^{3(p+q)+5}+u^{2(p+q)+6}\beta(Y_{p,q}^{\epsilon}),\\
\beta(A_4(g_4^{+}))=(u-1)u^{3(p+q)+6}\beta(Y_{p,q}^*)+(u-1)^2u^{3(p+q)+5}+u^{2(p+q)+6}\beta(Y_{p,q}^{c}),\\
\beta(A_4(g_4^{-}))=(u-1)u^{3(p+q)+6}\beta(Y_{p,q}^*)+3(u-1)^2u^{3(p+q)+5}+u^{2(p+q)+6}\beta(Y_{p,q}^{c}).\\
\end{array}
\end{displaymath}

We postpone the details of these computations for a moment. Note that 
$$\beta(A_4^{\epsilon}(h))=u^{3(p+q)+6}\beta(Y_{p,q}^*)+u^{2(p+q)+7}\beta(Y_{p,q}^{\epsilon})$$
for $h=h_7,h_8$ and $\epsilon \in \{1,-1\}$ by Lemma \ref{lem7}, and this polynomial is different from $\beta(A_4^{\epsilon}(g_4^{\pm}))$. Indeed the equality would imply that $\beta(Y_{p,q}^{\epsilon})$ is equal either to $(u-1)u^{p+q-1}$ or $3(u-1)u^{p+q-1}$. That can not happen thanks to Proposition \ref{rappel}. Therefore $g_4^{\pm}$ can not blow-Nash equivalent neither to $h_7$ nor to $h_8$.

Finally, let's prove that $g_4^{\pm}$ can not be blow-Nash equivalent to $h_6^{\pm}$. Consider $h_6^+$ first. The equality of the virtual Poincar\'e polynomial of $A_4(g_4^{\pm})$ and $A_4(h_6^+)$ would imply by Lemma \ref{lem7}.\ref{lem7.2} that
$$\alpha (u-1)^2u^{3(p+q)+5}+u^{2(p+q)+6}\beta(Y_{p,q}^c)=u^{2(p+q)+6}\beta(Y_{p+1,q}^c),$$
where $\alpha=1$ for $g_4^+$ and $\alpha=3$ for $g_4^-$. As a consequence $\beta(Y_{p+1,q}^c)-\beta(Y_{p,q}^c)$ should be equal to $\alpha (u-1)^2u^{p+q-1}.$
This can not happen since
\begin{displaymath}
\beta(Y_{p+1,q}^c)-\beta(Y_{p,q}^c)= \left \{
\begin{array}{ll}
(u-1)^2u^{p+q-1} +u^p(1-u) \textrm{~~if~~} p<q,\\
(u-1)^2u^{p+q-1} +u^{p-1}(u-1) \textrm{~~if~~} p \geq q\\
\end{array}\right.
\end{displaymath}
by Proposition \ref{rappel}.

Similarly $g_4^{\pm}$ and $h_6^-$ can not be blow-Nash equivalent because an equality of the virtual Poincar\'e polynomial of their spaces of arcs $A_4$ would imply the equality 
\begin{displaymath}
\alpha (u-1)^2u^{p+q-1}= \left \{
\begin{array}{ll}
(u-1)^2u^{p+q-1} +u^{q-1}(u-1) \textrm{~~if~~} p<q+1,\\
(u-1)^2u^{p+q-1} +u^{q}(1-u) \textrm{~~if~~} p \geq q+1.\\
\end{array}\right.
\end{displaymath}
This is forbidden by Proposition \ref{rappel}.

To complete the proof, it remains to check the value of the virtual Poincar\'e polynomial for $A_4(g_4^{\pm})$ and $A_4^{\epsilon}(g_4^{\pm})$ announced upper. The corresponding system of equations is (keeping the notations from section \ref{sectionD}):
\begin{displaymath}
\left\{ \begin{array}{lll}
Q(c_1)=0\\
a_1(b_1^2\pm a_1^2)+ \Phi(c_1,c_2)=0\\
2a_1b_1b_2+a_2b_1^2+3a_1^2a_2+Q(c_2)+\Phi(c_1,c_3)=\epsilon\\
\end{array} \right.
\end{displaymath}
Similarly to the proof of Lemma \ref{lem3}, we decompose this algebraic set in order to compute its virtual Poincar\'e polynomial. Note that the difference between $\beta(A_4^{\epsilon}(g_4^{+}))$ and $\beta(A_4^{\epsilon}(g_4^{-}))$ comes from the fact that, at the step $c_1^i=0$ for $i=1,\ldots,p+q$, the second equation becomes $a_1(b_1^2\pm a_1^2)=0$. If $a_1\neq 0$, then no solution exists for $g_4^+$ whereas for $g_4^-$ it gives a contribution of $2(u-1)$. Then if $a_1=0$, the third equation is the same in both cases.

The same method applies for $\beta(A_4(g_4^{\pm}))$.
\end{proof}
%%%%%%%%%%%%%%%%%%%%%%%%%%%%%%%%%%%%%%%%%%%%%%%%%%%%%%%%%%%%%%%%%%%%%%%%%%%%%%%%%%%%%%%%%%%%%
%%%%%%%%%%%%%%%%%%%%%%%%%%%%%%%%%%%%%%%%%%%%%%%%%%%%%%%%%%%%%%%%%%%%%%%%%%%%%%%%%%%%%%%%%%%%%
%%%%%%%%%%%%%%%%%%%%%%%%%%%%%%%%%%%%%%%%%%%%%%%%%%%%%%%%%%%%%%%%%%%%%%%%%%%%%%%%%%%%%%%%%%%%%
\section{Simple versus nonsimple singularities}\label{nonsimple}

Next step in the classification is to be able to distinguish between simple and nonsimple singularities. Namely, may a nonsimple Nash function germ belong to the same blow-Nash equivalence class as that of a simple one?

\begin{thm}\label{thm2} Let $f,g:(\mathbb R^{d},0) \rightarrow (\mathbb R,0)$ be Nash function germs. Assume $f$ is simple. If $f$ and $g$ are blow-Nash equivalent, then $g$ must be simple.
\end{thm}

For the proof of Theorem \ref{thm2}, we need first to reduce the possible types for the nonsimple germs before analysing more in details the possible spaces of arcs as in the proof of Theorem \ref{thm1}. 

\begin{cor} Let $f,g:(\mathbb R^{d},0) \rightarrow (\mathbb R,0)$ be Nash function germs. Assume $f$ is simple. Then  $f$ and $g$ are blow-Nash equivalent if and only if $f$ and $g$ are analytically equivalent.
\end{cor}

\begin{rmk}\label{rmkfinal} This result is maximal in the sense that it remains no longer true among nonsimple singularities. The family $f_t(x,y)=x^4 +2tx^2y^2+y^4+y^5$ defined on $(\mathbb R^2,0)$ is blow-Nash trivial along $t>0$ since its homogeneous part admits an isolated singularity \cite{fichou}. However the multiplicity is $9$ for $t\neq 1$ whereas for $t=1$ it becomes $11$.
\end{rmk}
%%%%%%%%%%%%%%%%%%%%%%%%%%%%%%%%%%%%%%%%%%%%%%%%%%%%%%%%%%%%%%%%%%%%%%%%%%%%%%%%%%%%%%%%%%%%%%%%%%%

In view of the proof of Theorem \ref{thm2}, we recall part of Arnold list of corank two nonsimple real isolated singularities \cite{Arnold}. We are interested in those singularities with nonzero 3-jet.

Let's denote $A_1=0$ and $A_k(x)=a_0+ a_1x+\cdots +a_{k-2}x^{k-2}$ for $k>1$ and $a_0,\ldots,a_{k-2}\in \mathbb R$. 
A nonsimple Nash function germ of corank two with an isolated singularity at the origin and a nonzero 3-jet is Nash equivalent to one of the following function germs:
\begin{displaymath}
\begin{array}{llllll}
J_{k,0}:  x_1^3 \pm bx_1^2x_2^{k}\pm x_2^{3k}+A_{k-1}(x_2)x_1x_2^{2k+1} \textrm{~~for $k>1$ and~~} 4b^3+27 \neq 0,\\
J_{k,i}:  x_1^3\pm x_1^2x_2^k+A_k(x_1)x_1^{3k+i} \textrm{~~for~~} k>1, i>0, a_0\neq 0, \\
E_{6k}:  x_1^3\pm x_2^{3k+1}+A_k(x_1)x_1x_2^{2k+1} \textrm{~~for~~} k> 1,\\
E_{6k+1}:  x_1^3+x_1x_2^{2k+1}+A_k(x_1)x_2^{3k+2} \textrm{~~for~~} k> 1,\\
E_{6k+2}:  x_1^3\pm x_2^{3k+2}+A_k(x_1)x_1x_2^{2k+2} \textrm{~~for~~} k> 1.\\
\end{array}
\end{displaymath}

Actually, the result holds in the real analytic setting, and the Nash case is a consequence of Nash Approximation Theorem (see section \ref{ADE}).

\begin{proof}[Proof of Theorem \ref{thm2}] 
To begin with, note that $g$ must have an isolated singularity similarly as $f$ because $f$ and $g$ are blow-Nash equivalent \cite{P}.

For the same reason they also share the same corank and index by Corollary \ref{invcorank}. Therefore we may suppose that $f$ and $g$ have the same quadratic part $Q$. Note that possible coranks for simple function germs are one or two. In case of corank one, then $g$ has an $A_k$-singularity and Theorem \ref{thm2} is proved. Otherwise, we may suppose that $f$ and $g$ are defined by
\begin{displaymath}
\begin{array}{ll}
f(x_1,x_2,y)=\widetilde f(x_1,x_2) + Q(y),\\
g(x_1,x_2,y)=\widetilde g(x_1,x_2) + Q(y),\\
\end{array}
\end{displaymath}
where $\widetilde f$ has a corank two simple singularity and $\widetilde g$ is a Nash germ with order at least 3. Note that $\widetilde f$ must be of type $D_k$ or $E_6,E_7,E_8$ as a corank two simple germ. Let's assume that $g$ is not simple.

We have proved in section \ref{ADE} that, for such a simple germ $f=\widetilde f+Q$ of corank two, the virtual Poincar\'e polynomial of the arc spaces $A_3^{\pm}(f)$ is not simply equal to the contribution coming from the quadratic part $(x_1,x_2,y)\mapsto Q(y)$, but an additional term does exist. Namely
$$\beta(A_3^{\pm}(f))=\beta(A_3^{\pm}(Q))+u^{2(p+q)}\beta(A_3^{\pm}(\widetilde f))$$
with $\beta(A_3^{\pm}(\widetilde f)) \neq 0$, where $(p,q)$ denotes the index of $f$. 

Therefore $\beta(A_3^{\pm}(g))$ must be different from $\beta(A_3^{\pm}(Q))$ because $f$ and $g$ share the same zeta functions. As a consequence $\widetilde g$ must have a nonzero 3-jet. According to Arnold list, the degree 3 homogeneous part $\widetilde g_3$ of $\widetilde g$ can be supposed to be equal to $\widetilde g_3(x_1,x_2)=x_1^3$, and $g$ must belong to one of the families of corank two nonsimple real isolated singularities with a nonzero 3-jet.

Moreover $g$ can not be blow-Nash equivalent to a $D_k$-singularity. Actually
$$\beta(A_3^{\pm}(x_1^3+Q(y,z)))=\beta(A_3^{\pm}(Q))+u^{2(p+q)+5}$$
is different from
$$\beta(A_3^{\pm}(G))=\beta(A_3^{\pm}(Q))+u^{2(p+q)+4}(u-1)$$
by Lemma \ref{lem4}, so $g$ can not be blow-Nash equivalent to a $D_k$-singularity for $k>4$. Moreover $\beta(A_3^{\epsilon}(x_1^3+Q(y,z)))$ is different from
$$\beta(A_3^{\epsilon}(x_1x_2^2-\epsilon x_1^3+Q(y,z)))=\beta(A_3^{\pm}(Q))+2u^{2(p+q)+5}$$
by Lemma \ref{lem5}, so $g$ can not either be blow-Nash equivalent to a $D_4$-singularity.

To complete the proof, it suffices to show that $g$ can not belong to the same blow-Nash class as that of a $E_6^{\pm},E_7$ or $E_8$ singularity. Assume that $f$ is of one of these types.

Analogously, the computations of Lemma \ref{lem7} prove that at least one of the virtual Poincar\'e polynomials of the arc spaces $A_4^{\pm}(f)$ or $ A_5^{\pm}(f)$ is not equal to the corresponding polynomial for the germ $h:(x_1,x_2,y)\mapsto x_1^3+Q(y)$, namely
$$\beta(A_4^{\pm}(h))=u^{3(p+q)+6} \beta(Y_{p,q}^*)+u^{2(p+q)+7} \beta(Y_{p,q}^{\pm})$$
and
$$\beta(A_5^{\pm}(h))=u^{4(p+q)+7} \beta(Y_{p,q}^*)+u^{3(p+q)+8} \beta(Y_{p,q}^*).$$

Therefore $\widetilde g$ must have a nonzero degree 4 or degree 5 homogeneous part since
$$\beta(A_4^{\pm}(f))=\beta(A_4^{\pm}(g)) \textrm{~~and~~} \beta(A_5^{\pm}(f))=\beta(A_5^{\pm}(g)).$$
As a consequence, the function germ $\widetilde g$ must belong to one of the families $J_{2,0},J_{3,0},J_{2,i},J_{3,i}$ for $i>0$. But, even for members of these families, the arc spaces $A_4^{\pm}(g), A_5^{\pm}(g)$ are the same as that of $h$. Indeed the composition with $\widetilde g$ of an arc $\gamma (t)=(a_1t+\ldots,b_1t+\ldots)$ is of the form
$$\widetilde g(\gamma (t))=a_1^3t^3+a_1^2(\cdots)t^4+a_1(\cdots)t^5+\cdots$$
so that the vanishing of $a_1$ implies the freedom of all other variables.

Consequently $g$ can not either be blow-Nash equivalent to a $E_6^{\pm},E_7$ or $E_8$ singularity.
\end{proof}
%%%%%%%%%%%%%%%%%%%%%%%%%%%%%%%%%%%%%%%%%%%%%%%%%%%%%%%%%%%%%%%%%%%%%%%%%%%%%%%%%%%%%%%%%%%%%%%%%%%%%%%%%%%%%%%%%

\end{document}